\documentclass[12pt]{article}
\usepackage[cp1251]{inputenc}
\usepackage[T2A]{fontenc}
\usepackage[russian,english]{babel}
\usepackage{amsfonts}
\usepackage{amsmath}
\usepackage{amssymb}
\usepackage{amscd}
\usepackage{graphicx}
\usepackage{color}
\newcounter{assumption}
\newtheorem{theorem}{Theorem}

\newtheorem{definition}{Definition}

\newtheorem{lemma}{Lemma}
\newtheorem{proposition}{Proposition}
\newtheorem{remark}{Remark}

\newenvironment{proof}[1][Proof]{\textbf{#1.} }{\ \rule{0.5em}{0.5em}}

\newcommand{\R}{{\mathbb R}}
\newcommand{\Z}{{\mathbb Z}}

\renewcommand{\dim}{\mbox{\rm dim}}

\title{Non-autonomous vector fields on $S^3$: \\simple dynamics
and wild separatrices embedding}
\author{V.Z. Grines$^{1}$, L.M. Lerman$^{1,2}$\\
\normalsize
$^1$Higher School of Economics, Russia\\
\normalsize
e-mail: vgrines@yandex.ru\\
\normalsize
$^2$Research and Educational Center ``Mathematics for Future Technologies'',\\
\normalsize
Lobachevsky National Research State University of Nizhny Novgorod,\\
\normalsize
e-mail: lermanl@mm.unn.ru \\
}
\date{}
\begin{document}
\maketitle

{\bf Abstract}. We construct new substantive examples of non-autonomous vector fields on
3-dimensional sphere having a simple dynamics but non-trivial topology. The construction is
based on two ideas: the theory of diffeomorpisms with wild separatrix embedding (Pixton,
Bonatti-Grines, etc.) and the construction of a non-autonomous suspension over a diffeomorpism
(Lerman-Vainshtein). As a result, we get periodic, almost periodic or even nonrecurrent vector
fields which have a finite number of special integral curves possessing exponential dichotomy
on $\R$ such that among them there is one saddle integral curve (with an exponential dichotomy
of the type (3,2)) having wildly embedded two-dimensional unstable separatrix and wildly
embedded three-dimensional stable manifold. All other integral curves tend, as
$t\to \pm \infty,$ to these special integral curves. Also we construct another vector
fields having $k\ge 2$ special saddle integral curves with tamely embedded two-dimensional
unstable separatrices forming mildly wild frames in the sense of Debrunner-Fox. In the case
of periodic vector fields, corresponding specific integral curves are periodic with the period
of the vector field, and they are almost periodic for the case of almost periodic vector fields.

\section{Introduction}

The aim of this paper is to present new substantive examples
of non-autonomous periodic, almost periodic or even nonrecurrent vector
fields on 3-dimensional sphere $S^3$. The main feature of constructed vector fields
is that they have some saddle integral curve $\gamma$ whose two-dimensional unstable and
three-dimensional stable invariant manifolds wildly embed in the extended phase
manifold $S^3\times\R$. Thus, this provides new invariants of an uniform equivalence of
non-autonomous vector fields (see definitions below). For the rest, from a dynamical viewpoint,
the vector field has a quite simple structure of its foliation into integral curves
(ICs, for briefness) in $S^3\times \R$.

It will be seen from the construction that the method allows one to get  non-autonomous
uniformly dissipative vector fields in $\R^3$ with a similar structure whose integral curves
enter to some cylindrical domain of the form $D^3 \times \R$, $D^3 \subset \R^3,$
with its boundary manifold $S^2\times \R$ being uniformly transversal
to integral curves.

To be more precise, let us recall the notion of a non-autonomous vector
field on a smooth ($C^\infty$) connected closed manifold $M$.
Let  $\mathcal V^{r}(M)$ be a Banach space of $C^r$-smooth, $r\ge 1,$ vector fields on
$M$ endowed with $C^r$-norm. By a \textit{$C^r$-smooth non-autonomous
vector field} on $M$ (NVF, for brevity) it is understood an uniformly continuous
bounded map $v: \mathbb R \to \mathcal V^{r}(M)$. We endow the set of non-autonomous
vector fields with the supremum norm of the related maps. As a particular case, one may
think on a periodic NVF, if the map $v$ is periodic: there exist a positive $T\in \mathbb
R$ such that $v(t+T)\equiv v(t)$ for all $t\in \mathbb R.$ If the map $v$
is almost periodic \cite{LZh,Cord}, they say on an almost periodic non-autonomous vector
field.

We recall {\em a solution} to the vector field $v$ be a $C^1$-differentiable map $x: I\to M$,
$I$ is an interval of $\mathbb R,$ such that for any $t\in I$ tangent vector
$x'(t)= Dx_t(1)\in TM$ coincides with the vector $v_t(x(t)).$
Here we identify in the standard way the tangent space
$T_t\mathbb R$ at the point $t\in \mathbb R$ with $\mathbb R$
itself by shifts in $\mathbb R$.

By the standard existence and uniqueness theorem, there is a unique solution through
any initial point $(x_0,\tau)\in M\times \R.$ On a closed manifold $M$ any solution of $v$
is extended on the whole $\mathbb R.$ The graph of the map $x$, that is the set
$\cup_t(x(t),t)\subset M\times \mathbb R$, is the integral curve of the solution $x$.
Thus, every non-autonomous vector field $v$ generates
a foliation $\mathcal L_v$ of manifold $M\times \mathbb R$ into ICs of $v$.
An example of such foliation for the case $M = I$ is plotted on Fig.~\ref{1dim}.

Following \cite{LSh} we call two non-autonomous vector fields $v_1, v_2$ on $M$ to be
{\em uniformly equivalent}, if there is an equimorphism (see subsection \ref{Addendum}
for definition) $\Phi:M\times \mathbb R \to M\times \mathbb R$ that transforms
$\mathcal L_{v_1}$ to $\mathcal L_{v_2}$ preserving orientation in $\R.$
Here we consider manifold $M\times \mathbb R$ with the uniform structure of the
direct product of the unique uniform structure on $M$ given by the topology of the manifold
$M$ (it is a compact manifold) and the standard uniform structure on $\R$ invariant w.r.t.
shifts on the Abelian group (see subsection \ref{Addendum} for definitions of uniform
structure and equimorphism).

\begin{figure}[h]
\centerline{\includegraphics[height=10cm]{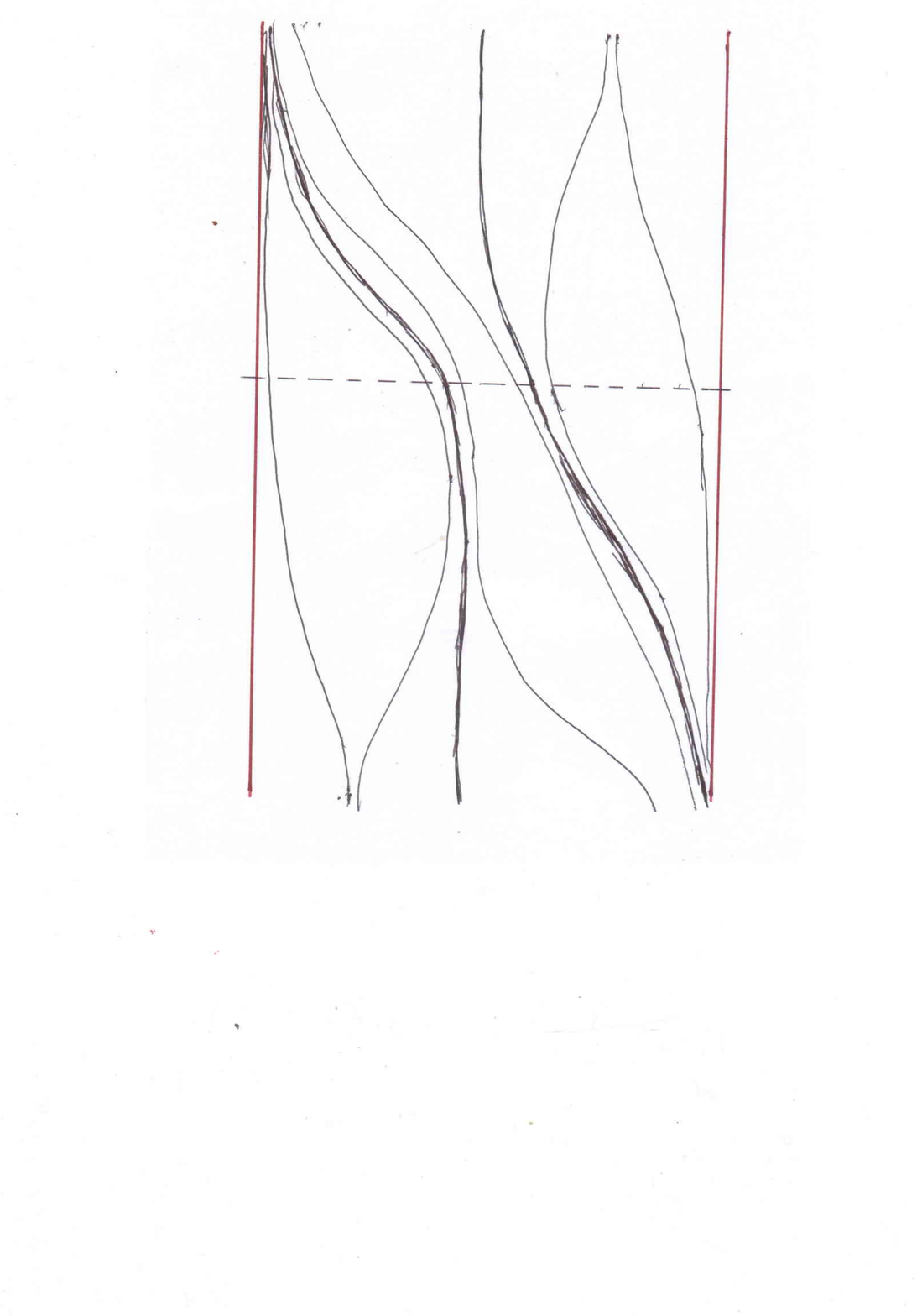}}
\caption{\small A foliation into IC on the segment, $M=I.$}
\end{figure}
\label{1dim}

NVFs, we construct below, fall into the class of gradient-like non-autonomous
NVFs singled out in \cite{LSh,Ler_dis}. They satisfy several restrictions on the structure of
their foliation $\mathcal L_v$, one of them is the claim of an exponential dichotomy
for any of their solutions on both semi-axes $\R_+$ and $\R_-$ (types of dichotomy
may distinct) \cite{MS}. This assumption allows one to get invariant stable manifolds
and also unstable manifolds. Thus, the whole extended phase manifold $M\times \R$ is partitioned
into smooth stable manifolds, another partition is generated by unstable manifolds. One more
assumption is the finiteness of both partitions (though they can be completely different).
We shall not deep into details of these restrictions, since in the examples we construct
these restrictions are given explicitly.

The NVF on $S^3$, we shall construct, has four special ICs possessing an exponential
dichotomy on the whole $\R.$ One such IC $\gamma_\alpha$ is exponentially unstable on
$\R$, there is also the only IC $\gamma_\sigma$ of a saddle type, it possesses an exponential
dichotomy on $\R$ of the type (3,2), that is, such IC has 3-dimensional stable manifold
and 2-dimensional  unstable manifold on $\R$. Finally, this NVF has two exponentially stable
on $\R$ ICs $\gamma_{\omega_1},$ $\gamma_{\omega_2}$ whose stable manifolds are 4-dimensional
ones (the dichotomy of the type (1,4)).

The sense of the term ``wild embedding'' we use below for some embedded submanifold is
the following. Take any section $t= t_0$ in $S^3\times \R.$ Then the intersection of the
wildly embedded 2-dimensional unstable manifold with the section is a 1-dimensional ray
in $S^3\times \{t_0\}$ wildly embedded in the topological sense \cite{AF,Pixton1977,BG}
(see below).
The closure of this ray is a point being the trace of the exponentially stable IC
$\gamma_{\omega_2}$. All ICs in this unstable manifold tend to $\gamma_\sigma$
as $t\to -\infty$ and to $\gamma_{\omega_2}$ as $t\to \infty$. Also, the trace
of 3-dimensional stable manifold of $\gamma_\sigma$ on $M_{t_0}$ is an embedded
$\R^2$ and its closure in $M_{t_0}$ is an embedded sphere being wild at one point \cite{}.

\begin{figure}[h]
\centerline{\includegraphics[height=8cm]{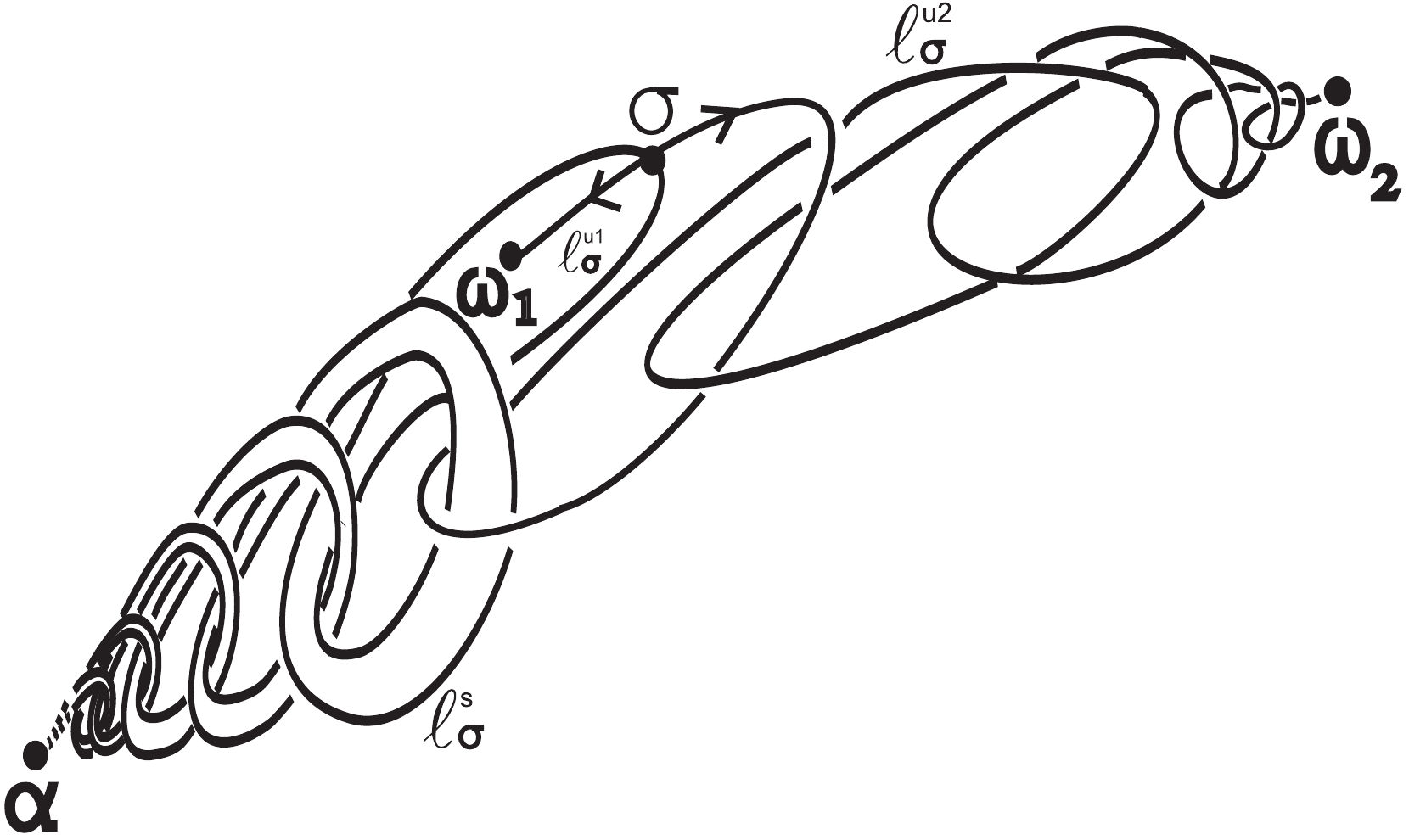}}\caption{\small Phase portrait
of diffeomorphism on $S^3$ with wildly embedded separatices.}
\label{p}
\end{figure}

The construction of such non-autonomous vector fields exploits two ideas. One
belongs to Pixton \cite{Pixton1977} and further was developed by Bonatti, Grines, Pochinka
and others \cite{BG,BGMP,GMP}. In \cite{Pixton1977} a simple 3-dimensional Morse-Smale
diffeomorphism on $S^3$ was constructed that has, as its non-wandering set, only
four hyperbolic fixed points: one source $\alpha$, one saddle $\sigma$ of $(2,1)$ type and
two sinks $\omega_1$, $\omega_2$. Moreover, the closure of the stable manifold $w^s_{\sigma}$
of $\sigma$ is homeomorphic to the sphere smoothly  embedded in each point except
the point $\alpha$ where it is wildly embedded   (see, on the figure \ref{p}).
The point $\sigma$ divides unstable manifold $w^u_\sigma$ into two separatrices
$l^u_{1},l^u_{2}$. The closure of one of these separatrices ($l^u_{2}$ on
the figure \ref{p}) is a simple arc wildly embedded  in the point $\omega_2$  but
the closure of  another one ($l^u_{1}$) is tamely embedded arc in each point.
The separatrix $l^u_{2}$ tends to the sink $\omega_2$ in such a way that the fundamental
domain near this sink after identifying points on the boundary
of this domain contains the image of the separtrix which makes up a nontrivial knot
in the related factor space (being the manifold $S^2\times S^1$). This will be explained
in more details in section \ref{wild-sep}. The non-triviality of this knot implies the wild
embedding of the separatrix for the diffeomorphism, moreover,  if knots for two different diffeomorphisms are non-homeomorphic
in $S^2\times S^1$, then they are not topologically conjugated. This  fact leads to the existence of infinite number of topologically non-conjugated Morse-Smale
diffeomorpisms of considered types   \cite{BG}.

The second idea is borrowed from \cite{LV}. It uses the construction of a
so-called non-autonomous suspension over a diffeomorphism introduced in \cite{LV0} and
developed further in \cite{LV}. Recall this construction. Suppose $f: M\to M$ be some
diffeomorphism of a smooth ($C^\infty$) closed manifold $M$. To avoid a discussion on the
class of smoothness for the suspension, we assume $f$ to be $C^\infty$-smooth. Its (usual)
suspension is a smooth closed manifold $M_f$ of dimension $\dim M +1$ with a flow
defined as follows. Let us identify in the cylinder $M\times I,$ $I=[0,1],$ points $(x,1)$
and $(f(x),0)$. It is more convenient to consider the manifold $M\times \R$ with an action
$F$ of the group $\Z$ by the rule: for $m\in \Z$ the related $F^m$ acts
as $F^m(x,s) = (f^m(x),s-m)$. This action is free and discrete (any orbit of the action has
not accumulation points). Thus, the factor-manifold $M_f = (M\times \R)/F$ is
a smooth manifold being a smooth bundle over the circle $S^1$, $p:M_f \to S^1$, with the leaf
$M$. A vector field on $M_f$ is generated by the constant vector field $V=(0,1)$ on $M\times
\R$ (its orbits are straight-lines $(x,t)$, $t\in \R$). After factorizing one
gets a smooth vector field $v_f$ on $M_f$ with the global cross-section, as such one can
choose any $M_\theta= p^{-1}(\theta),$ $\theta \in S^1$. The Poincar\'e map defined on
this cross-section is conjugated to $f.$ This construction allows one to get vector
fields with the dynamics similar to that as for iterations of the mapping $f$, see \cite{Sm}.

To proceed let us consider the covering manifold $\tilde{M}_f$ for $M_f$ generated by the
standard covering $\R \to S^1$, $t\to \exp[2\pi i t]$ that gives a
commutative diagram
$$
\begin{CD}
\tilde{M}_f@>\widetilde{\exp}>> M_f\\
@VV\tilde{p}V            @VVpV\\
\R@>\exp>>               S^1
\end{CD}
$$
Topologically, $\tilde{M}_f$ is homeomorphic to $M\times \R$, since $\R$
is contractible. The manifold $M_f$ is compact, hence it has the unique uniform structure
compatible with the topology \cite{Kelley}. The  uniform  structure  in  $\tilde{M}_f$
is defined by lifting the uniform structure in $M_f$ by means of the map
$\widetilde{\exp}$. This is more easily understood, if we endow $M_f$ with a smooth Riemannian
metrics and lift this metrics to $\tilde{M}_f$ by the covering map $\widetilde{\exp}$.
Since $\widetilde{\exp}$ is a local diffeomorphism, we get a Riemannian metrics on $\tilde{M}_f$
such that $\widetilde{\exp}$ is the local isometry.
The foliation in $M_f$ into orbits of the vector field $v_f$ is lifted as some foliation
$\mathcal L_{v_f}$ into infinite curves in $\tilde{M}_f$. This foliation is homeomorphic
to the foliation of $M\times \R$ into straight-lines $(x,t)$, $t\in\R$, since $\tilde{M}_f$
is homeomorphic to $M\times \R$, but generically the foliation in $\tilde{M}_f$ is not
equimorphic to the foliation into straight-lines.
Moreover, even the manifold $\tilde{M}_f$ itself with its uniform structure lifted from $M_f$
is not always equimorphic to $M\times \R$. For instance, it is the case for $M=T^2 =
\R^2/\Z^2$ with $f$ being an Anosov diffeomorphism (see details in
\cite{LV}). Next proposition which is Corollary 4.1 from \cite{LV}
will be useful for us.
\begin{proposition}
If $\tilde{M}_f$ is equimorphic to $M\times \R$, then there is such $n\in \mathbb
Z$ that $f^n$ is homotopic to $id_M.$
\end{proposition}

Let a diffeomorphism $f: M\to M$ on a smooth ($C^\infty$) closed manifold $M$ be given.
An important question here is if there exists a non-autonomous vector field $v$ on $M$ such
that its foliation $\mathcal L_v$ into ICs in $M \times \R$ (with its uniform structure
of the direct product) is equimorphic to the foliation $\mathcal L_{v_f}$ into infinite curves
generated by vector field $v_f$ in $\tilde{M}_f$?
It is evident the first condition this to be true is an equimorphness of uniform spaces
$M\times \R$ and $\tilde{M}_f$. This gives a meaning to the definition introduced
in \cite{LV}
\begin{definition}
A diffeomorphism $f: M\to M$ is reproduced by a non-autonomous vector field
$v$ on $M$ $($or, equivalently, $v$ reproduces the structure of $f$$)$,
if foliations $\mathcal L_v$  in $M\times \R$ and $\mathcal L_{v_f}$ in $\tilde{M}_f$
are equimorphic. In particular, this implies the uniform spaces $M\times \R$ and $\tilde{M}_f$
be equimorphic.
\end{definition}
\begin{remark}
It follows from Proposition 2.5 (item a) in $\cite{LV}$ that diffeomorphisms $f$ and $f^n$ are
reproduced simultaneously for any $n\in\Z$.
\end{remark}\label{remark-simult}
\begin{remark}
As is known, for a given diffeomorphism $f: M\to M$ it is possible that its suspension $M_f$
is not diffeomorphic to the direct product $M\times S^1$ but for some its
iteration $f^n$ the related suspension $M_{f^n}$ is diffeomorphic to $M\times S^1$. In fact, the manifold $M_{f^n}$ is the k-fold covering of $M_{f}$. As a simplest
example of such situation, one can take $M=S^1$ with a coordinate $\varphi\;(\mbox{\rm\,mod}\;
2\pi)$ and a diffeomorpism $f(\varphi) = -\varphi (\mbox{\rm\,mod}\;2\pi).$ Then $S^1_f$ is
the Klein bottle but $f^2(\varphi)= \varphi$, hence $S^1_{f^2} = S^1\times S^1 = \mathbb T^2$
(that is, a 2-dimensional torus being a 2-fold covering of the Klein bottle).
Thus, according to  remark \ref{remark-simult} the manifold  $\tilde{S^1_f}$ is
equimorphic to $S^1\times \R.$
\end{remark}\label{remark-tor-bottle}

Denote $\pi_M: M\times \R \to M$ the standard projection onto the first co-factor.
For any non-autonomous vector field $v$ on the manifold $M$ the map $\Phi^t_0: M_0 \to
M_t$ from the section $M_0 = M\times \{0\}$ to the manifold $M_t = M\times \{t\}$,
generated by solutions of $v$ with initial points on $M_0$, $\pi_M(M_t)= M,$ is diffeotopic to
the identity map $id_M$ for all $t\in \R$. In particular, if $v$ is a periodic vector field
on manifold $M$, then its Poincar\'e map over the period is diffeotopic to
identity map $id_M.$

Our first result is theorem \ref{T1} below  which gives sufficient condition for
diffeomorphism $f$ to be reproduced by a flow  generated by a non-autonomous periodic
vector field $v$.

First we formulate an obvious lemma.
\begin{lemma}\label{sm}
If $f$ is diffeotopic to $id_M$, then there is a diffeotopy $F_t: M\to M,$
$t\in [0,1],$ joining $id_M$ and $f$ such that diffeomorpisms $F_t$ depend smoothly
on $t$ and for some $\varepsilon > 0$ small enough one gets $F_t \equiv
id_M,$ as $t\in [0,\varepsilon],$ and $f_t \equiv f$ as $t\in [1-\varepsilon, 1].$
\end{lemma}

\begin{theorem}
Suppose for some $n\in \mathbb N$ diffeomorphism $f^n: M\to M$ be diffeotopic to
the identity map $id_M$. Then
\begin{enumerate}
\item $M_{f}$ is fiber-wisely\footnote{The term ``fiber-wisely'' means the existence
of a diffeomorphism $\Psi: M_{f} \to M\times S^1$ acting as $(x,s)\to (\psi(x,s),s).$}
diffeomorphic to $M\times S^1$;
\item there is a periodic vector field $v$ on $M$ such that $v$ reproduces
the structure of $f$.
\end{enumerate}
\end{theorem}\label{T1}
\proof
To ease the exposition we assume that $f$ itself is diffeotopic to $id_M$.
At the first step we shall construct a 1-periodic vector field $v_t$ on $M$ such that
the vector field on $M\times S^1$ given as $(v_t,1)$ is diffeomorphic
to the vector field of the suspension of $f$ on $M_f$. To this end, we need to endow
the manifold $M_f$ with the structure of the direct product in the explicit form.
This means that one needs to define two foliations of $M_f$. One of them is given
by leaves of the bundle over $S^1$, these leaves are diffeomorphic to $M$. The second
foliation into closed curves is defined as follows.
Suppose $F_t: M\to M,$ $t\in [0,1],$ be a diffeotopy $F_t: M\to M$ joining $id_M$
and $f$, that is, $F_t$ are diffeomorphisms, $F_0 = id_M$ and $F_1 = f.$ We assume,
by Lemma \ref{sm}, that $F_t = id_M$ for $t\in [0,\varepsilon]$ and $F_t = f$ for
$t\in [1-\varepsilon,1].$
If $p: M_f \to S^1$ is the bundle map, then for any point $x\in p^{-1}(0)$ we define
the curve through point $(x,0)\in M\times [0,1]$ given as $(F^{-1}_t(x),t)$ for
$t\in [0,1]$, and then we apply the factor map using the identification $(x,t)=(f(x),t-1).$
The curve in $M\times \R$ with the initial point $(x,0)$ has the
extreme point $(F^{-1}_1(x),1)=$ $(f^{-1}(x),1)$ at $t=1$. After identifying this
point becomes $(f\circ f^{-1}(x), 0)=(x,0).$ Thus, all curves constructed
in the manifold $M_f$ are closed and we get the
homeomorphism $h:M_f\to M\times S^1.$ The map $h$ is defined as follows.
Take any point $a\in M_f$ and denote $l_a$ a closed curve through the point $a$ of
the second foliation constructed. Define a map $p_1: M_f \to M_0$, $M_0 = p^{-1}(0),$
by the rule $p_1(a)= l_a\cap M_0$. We get a homeomorphism $h: a\to (p_1(a),p(a))$.

Generically, $h$ is a homeomorphism, since the dependence of $F_t$ on $t$ can be
only continuous. But we need a diffeomorpism between bundles $M_f$ and $M\times S^1$
in order orbits of the suspended flow would be transformed to smooth
curves in $M\times S^1$. Lemma \ref{sm} above guarantees that if
$f:M\to M$ is diffeotopic to $id_M,$ then a diffeotopy $F_t$ exists joining $id_M$
and $f$ such that curves constructed above give a smooth foliation, that
is, curves are smooth and their dependence of the points is smooth, so the map $p_2$
is smooth. This proves the first item of the theorem.

Now we construct a periodic vector field $v$ on $M$ such that its
foliation $\mathcal L_v$ into ICs is uniformly diffeomorphic to the
foliation $\mathcal L_{v_f}$ into infinite curves in $\tilde{M}_f.$
Diffeomorphism $h: M_f \to M\times S^1$, defined above, allows us to
identify $M_f$ with $M\times S^1.$ Thus, the suspended vector field is
given as $(V(x,t),n(x,t))$ with $V,n$ being 1-periodic in $t$ and $n > 0.$ Its
flow therefore has a cross-section, for instance, as such one can take the section
$t=0$. The Poincar\'e map $g:M_0 \to M_0$ on this cross-section is evidently conjugated to
$f$. Hence, we can consider instead of diffeomorphism $f$ on $M$ the
diffeomorphism $g$. Since $f,g$ are conjugated, their non-autonomous
suspensions are equimorphic along with the related foliations.
Compactness of $M$ and $S^1$ implies that $n$ is strictly positive.
We can define a periodic vector field on $M$ as $v(x,t)= V(x,t)/n(x,t)$.
Integral curves in $M\times \R$ of this periodic vector field coincide with orbits
the vector field $(V(x,t),n(x,t))$ since they are obtained by the change of time
being uniformly bounded from above and below. So, the item 2 has also been
proved.
$\blacksquare$

\section{Diffeomorphisms with wildly embedded separatricies}
\label{wild-sep}

This section contains some definitions and results which are contained in
the book \cite{GMP}, we present them here for  of the reader convenience.
\subsection{Wild embedding}

\begin{definition}\label{diko} A topological embedding $\lambda: X\to Y$ of
an $m$-dimensional manifold $X$ into a $n$-dimensional manifold $Y$ ($m\leq n$)
is said to be {\it locally flat at the point $\lambda(x)\in Y$}, if there is
a chart $(U,\psi),$ $\lambda(x)\in U,$ $\psi: U\to \R^n$, in the manifold $Y$ such that
$\psi(\lambda(X)\cap U) = D^m\subset\mathbb R^m$, here $\mathbb
R^m\subset\mathbb R^n$ is the set of points for which the last $n-m$ coordinates equal
to zero or $\psi(\lambda(X)\cap U)=\mathbb R^m_+$ $($$\mathbb
R^m_+\subset\mathbb R^m$ is the set of points with non-negative last coordinate$)$.
\end{definition}
Now, an embedding $\lambda$ is said to be {\em tame} and the manifold $X$ is said to be
{\em tamely embedded}, if $\lambda$ is locally flat at every point $\lambda(x)\in Y$.
Otherwise, the embedding $\lambda$ is said to be {\em wild} and the manifold $X$ is
said to be {\em wildly embedded}. If the embedding $\lambda$ is not locally flat at
the point $\lambda(x)$, this point is said to be a {\em point of
wildness}.

It is worth remarking that the definition of a tamely embedded manifold coincides with
the definition of a topological submanifold.

Every topological embedding into the space $\mathbb R^2$ (respectively, $S^2$) is tame.
In the space $\mathbb R^3$ (respectively, $S^3$) there are wild arcs and wild 2-spheres.
As an example of a wild arc, we recall the construction by Artin and Fox \cite{AF}.
The related arc is smooth everywhere except for its boundary point.

\begin{figure}[h]
\centerline{\includegraphics[height=4cm]{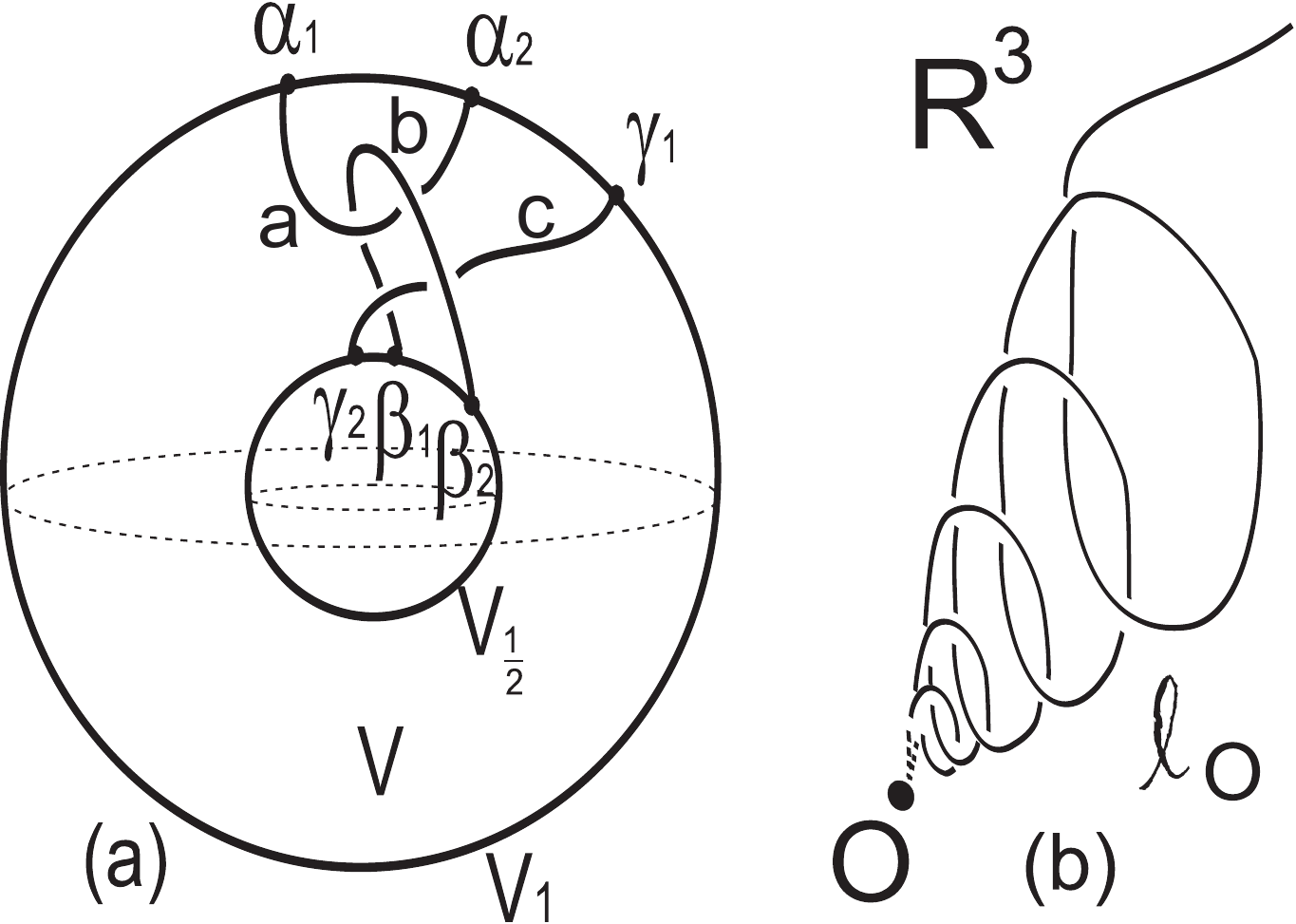}}\caption{{\small Constructions of
wild curves in $\mathbb R^3$}}
\label{r7}
\end{figure}

Consider a linear contraction $\phi:\R^3\to\R^3$ defined in spherical coordinates
$(\rho,\varphi,\theta)$ as $\phi(\rho,\varphi,\theta)= (\frac12\rho,\varphi,\theta)$,
and denote $L\subset \R^3$ a spherical layer defined by inequalities $\frac12\leq\rho\leq 1$.
Its boundary spheres are $V_{\frac12}=\{(\rho,\varphi,\theta)|\rho=\frac12\}$ and
$V_1=\{(\rho,\varphi,\theta)|\rho=1\}$. Let $a,b,c\subset L$ be pairwise disjoint simple
arcs with their respective boundary points $\alpha_1,\alpha_2;\beta_1,\beta_2;\gamma_1,\gamma_2$
(see Figure \ref{r7} (a)) such that
\begin{enumerate}
\item $\alpha_1,\alpha_2,\gamma_1\subset V_1$;
$\beta_1,\beta_2,\gamma_2\subset V_{\frac{1}{2}}$;
\item $\phi(\alpha_1)=\gamma_2$, $\phi(\alpha_2)=\beta_1$, $\phi(\gamma_1)=\beta_2$.
\end{enumerate}

 Let us  choose  arcs $a, b, c$ in     such a way that the arc $\ell_O\subset\mathbb R^3$
defined as $\ell_O=\bigcup\limits_{k\in\mathbb Z}\phi^k(a\cup b\cup c)\cup O$ (see Figure
\ref{r7} (b)), is  smooth at every point except $O$. Artin and Fox proved that $\ell_O$ is
wildly embedded into $\mathbb R^3$ and the point  $O$ is point of wildness. This fact also
follows from the criterion below proved in \cite{HGP}.

\begin{proposition}  Let $\ell$ be a compact arc in $\R^3$ which is smooth everywhere except
its boundary point $O$. Then $\ell$ is locally flat at $O$, iff for every
$\varepsilon$-ball $B_\varepsilon(O)$ centered at $O$ there is a subset
$U\subset B_\varepsilon (O)$ diffeomorphic to the closed 3-ball such that $O$ is an interior
point of $U$, and the intersection $\partial U\cap\ell$ is the only point.
\label{tame-arc}
\end{proposition}

Now we consider a standard sphere $S^3 = \{x = (x_1,x_2,x_3,x_4)\in \mathbb {R}^4|\;
x^2_1+x^2_2+ x^2_3+x^2_4=1\}$. The point $N(0,0,0,1) \in S^3$ (respectively,
$S(0,0,0, -1)\in S^3$) will be called the north (respectively, the south) pole.

\begin{figure}
\centerline{\includegraphics[height=4cm]{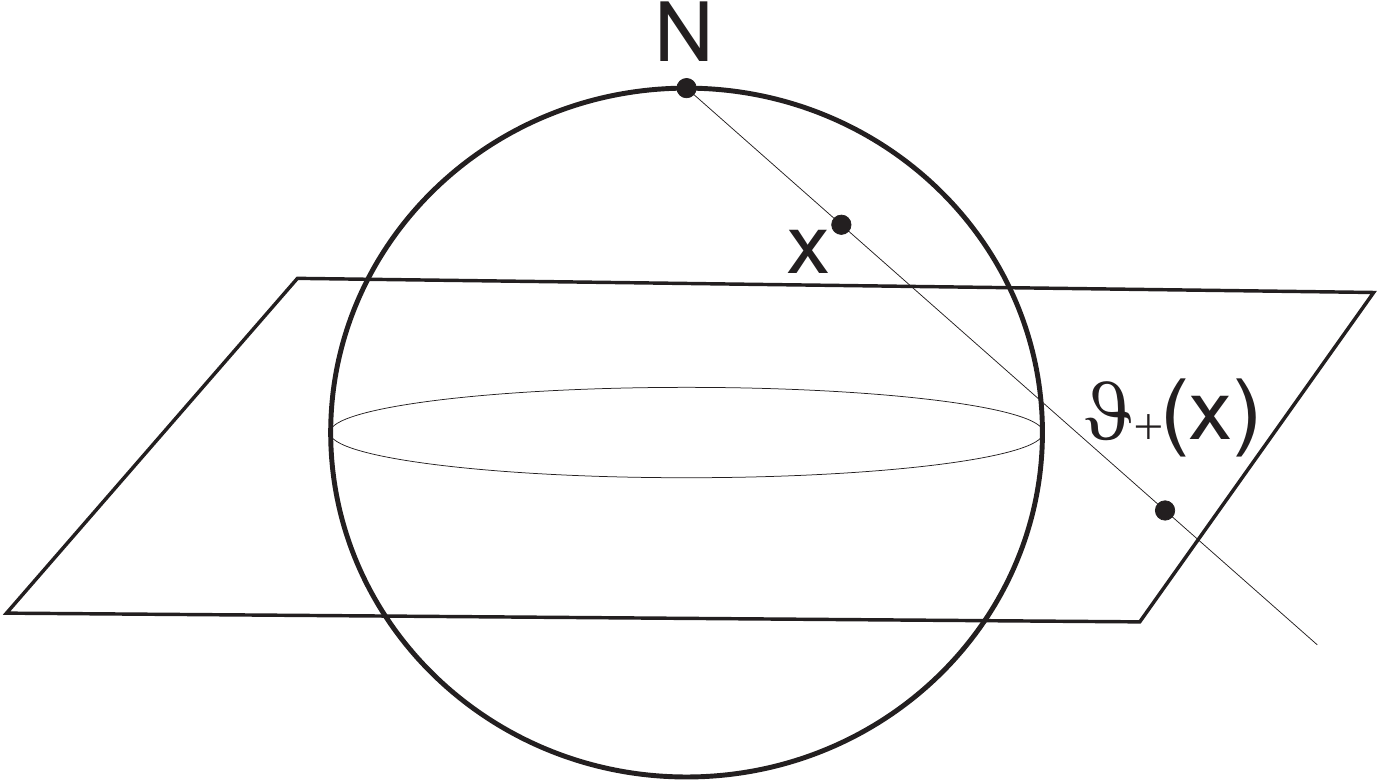}}\caption{{\small The stereographic
projection}}
\label{r9}
\end{figure}

For each point $x\in S^3\setminus\{N\}$ there is the unique straight line
in $\mathbb {R}^{4}$ containing $N$ and $x$. This line cuts the plane
$(x_1,x_2,x_3,0)$ at exactly one point
$\vartheta(x)$,  $\vartheta(x_1,x_2,x_3,x_{4})=\frac{x_1}{1-x_4},
\frac{x_2}{1-x_4},\frac{x_3}{1-x_4}$.
 The stereographic projection of the point $x$ is defined as the point $\vartheta(x)$.
The stereographic projection is a diffeomorphism of $S^3\setminus\{N\}$ to
$\mathbb R^3$ (see Figure \ref{r9}, where it is shown the stereographic projection
from $S^2\setminus\{N\}$ to $\mathbb R^2$).

\begin{figure}[h]
\centerline{\includegraphics[height=4cm]{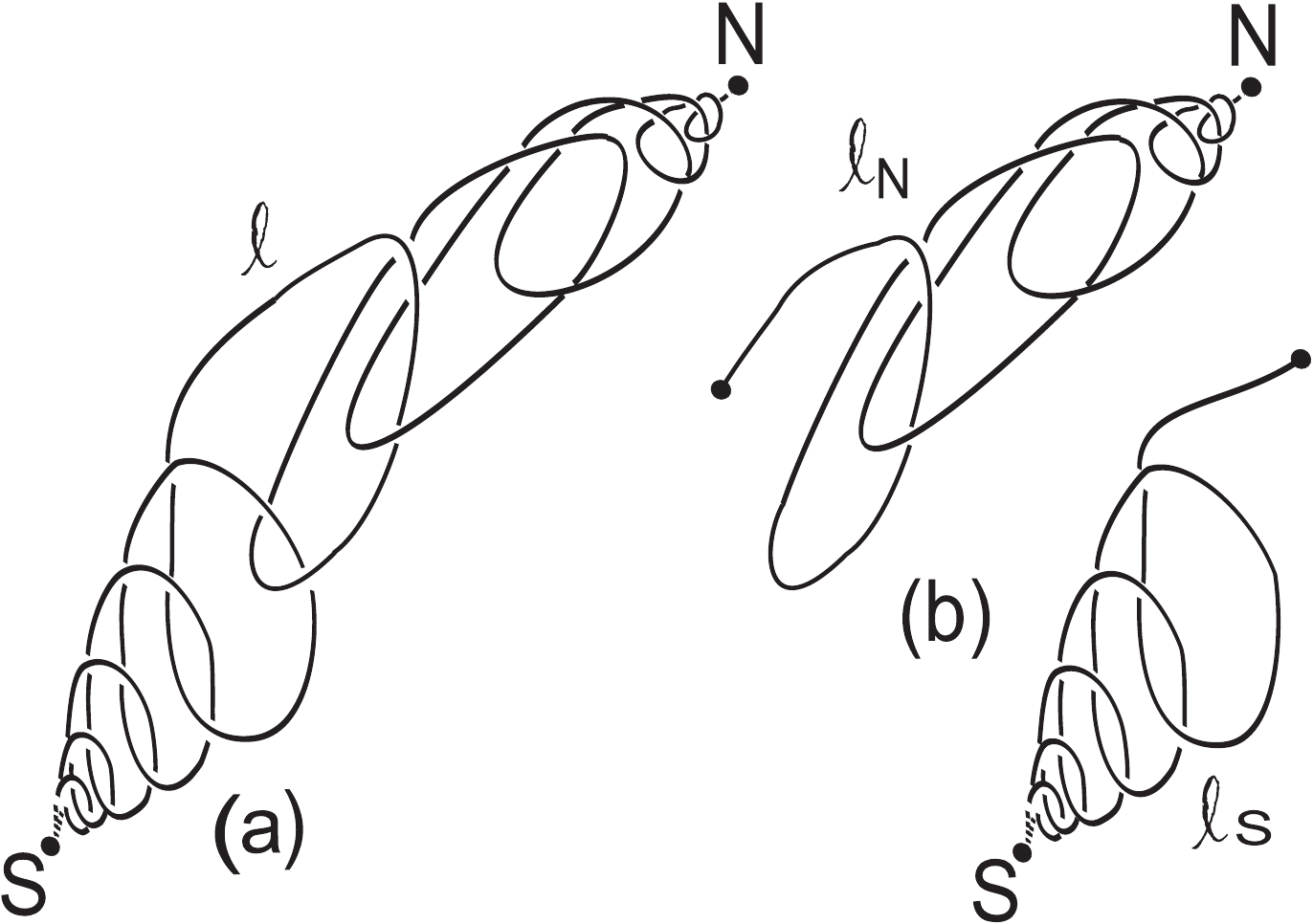}}\caption{{\small Constructions of
wild curves in $\mathbb S^3$}}
\label{r8}
\end{figure}

Let $\ell=\vartheta_+^{-1}(\ell_O)\cup S$ (see Figure \ref{r8} (a)), then the arc
$\ell_N$ ($\ell_S$) in Figure \ref{r8} (b) is a sub-arc of the arc $\ell$ from the point
$\vartheta_+^{-1}(\alpha_1)$ to the point $N$ (from the point $\vartheta_+^{-1}
(\alpha_1)$ to the point $S$). The arc $\ell_N$ ($\ell_S$) is wildly embedded into
$\mathbb S^3$.

Now we can inflate the arcs on Figure \ref{r8} (a), (b) to get closed 3-balls whose
boundaries are 2-spheres wildly embedded to $\mathbb S^3$ and whose wild points are
the poles.

\subsection{Diffeomorphisms of Pixton type on $S^3$}
\label{Pixton}

Let $V$ be a smooth closed orientable 3-manifold whose fundamental group admits
a nontrivial homomorphism $\eta_{_V}:\pi_1(V)\to\mathbb Z$. Denote by
$(V,\eta_{_V})$ the manifold $V$ equipped with the homomorphism
$\eta_{_V}$.

\begin{definition}\label{defi1} Manifolds $(V,\eta_{_V})$ and
$(V',\eta_{_{V'}})$ are said to be {\em equivalent}, if
there is a homeomorphism $\varphi: V\to V'$ such that
$\eta_{_{V'}}\varphi_* = \eta_{_V}$.
\end{definition}

\begin{definition}\label{defi2} Two smooth submanifolds $a\subset V$
and $a'\subset V'$ are said to be equivalent, if
there is a homeomorphism $\varphi: V\to V'$ such that
$\eta_{_{V'}}\varphi_*=\eta_{_V}$ and $\varphi(a) = a'$.
\end{definition}

\begin{definition}\label{defi3} A smooth submanifold $a\subset V$
is said to be $\eta_{_{V}}$-essential, if $\eta_{_V}(i_{a*}
(\pi_1(a)))\neq 0$, where $i_{a}: a\to V$ is the inclusion map.
\end{definition}

Let us illustrate these definitions for the manifold $S^2\times S^1$.
We represent the manifold $S^2\times S^1$ as the orbit space of the homothety
$a^s\mapsto a^s(x)=0.5x$ ($x = (x_1,x_2,x_3)$), $({\mathbb R}^3\setminus \{O\})/{a^s}$.
It is  easy to check that the natural  projection $p:
{\mathbb {R}}^3\setminus O\to S^2\times S^1$ is the covering map, it
induces the epimorphism $\eta_{_{S^2\times S^1}}:
\pi_1(S^2\times S^1)\to\mathbb{Z}$.\footnote{Take the homotopy
class $[c]\in\pi_1(S^2\times S^1)$ of a loop $c\colon \mathbb R/\mathbb Z
\to S^2\times S^1)$. Then $c\colon [0,1]\to S^2\times S^1$ lifts to a curve
$c\colon [0,1]\to \mathbb R^3\setminus \{O\}$ joining a point $x$ with a point
$(a^{s})^n(x)$ for some  $n\in \mathbb Z$, where $n$ is independent of the lift.
So, we define $\eta^s_{_{S^2\times S^1}}([c])=n$.}

\begin{figure}[h]
\centerline{\includegraphics[height=4cm]{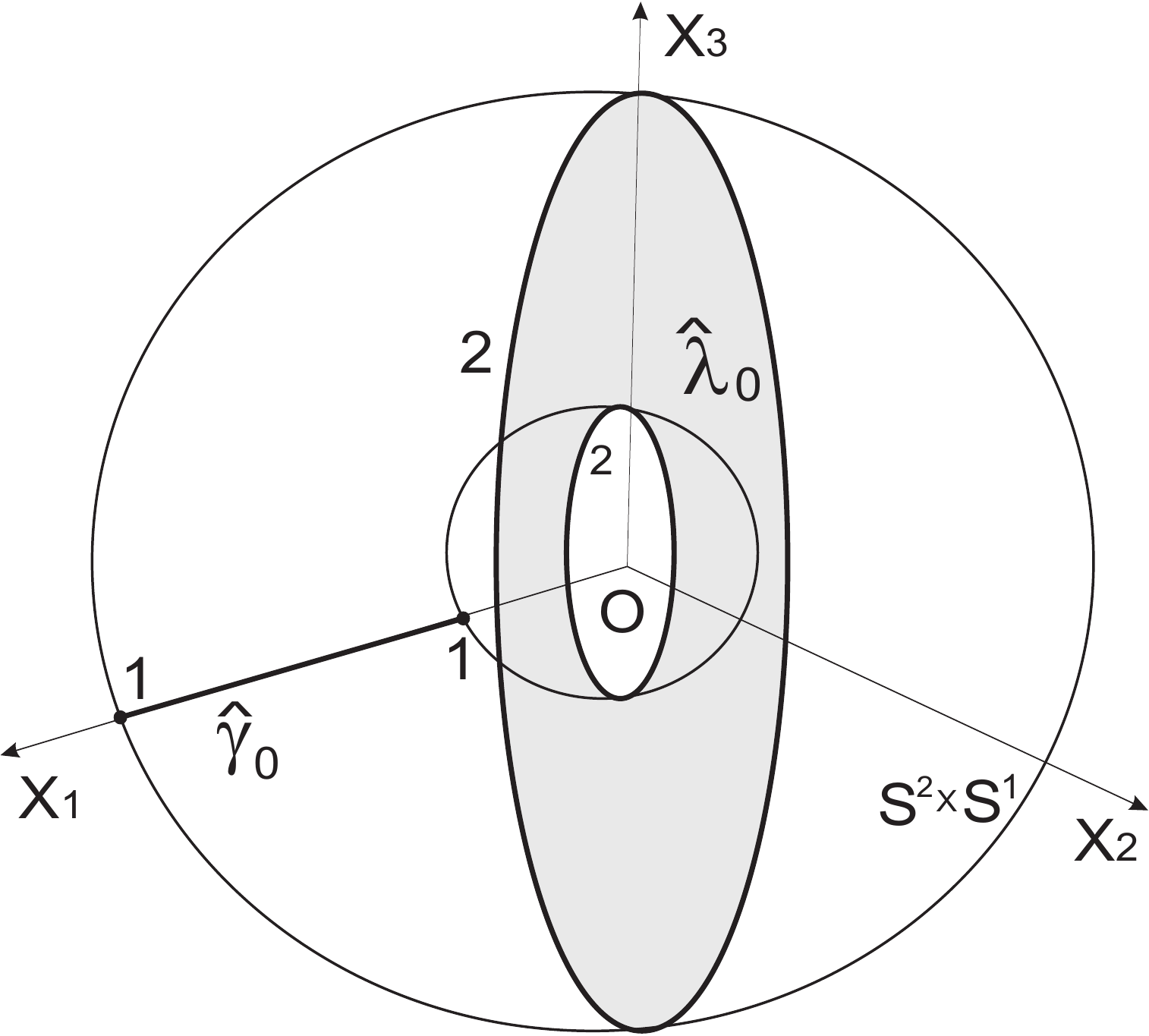}}\caption{{\small Construction of
an essential knot and torus embedded in $S^2\times S^1$}}
\label{r3}
\end{figure}

Denote $\gamma_0=p(Ox_1^+),\lambda_0 = p(Ox_2x_3)$, where $Ox_1^+$ is positive
semi-axis and $Ox_2x_3$ is coordinate plane $x_1=0$. On Fig.\ref{r3} it
is shown the spherical layer bounded by spheres of radii $1$ and $0.5$.  If we identify
points which lie on the boundary of the spherical layer and belong to the same
ray through $O$, we get the manifold $S^2\times S^1$. Moreover, if we identify the
extreme points of the segment with the same numbers ($1$), we get the knot
$\hat\gamma_0$ and if we identify extreme points lying on the same ray that belongs to
circles with the same numbers ($2$) and bounding $2$-annulus, we get the torus
$\hat\lambda_0$ ($\hat\gamma_0$ and $\hat\lambda_0$ are embedded to
$S^2\times S^1$).

It is easy to check that $\hat\gamma_0$ (respectively, $\hat\lambda_0$) is a
$\eta{_{\mathbb{S}^2\times\mathbb
{S}^1}}$-essential knot (respectively, torus) in the manifold $(S^2\times S^1,
\eta{_{\mathbb{S}^2\times\mathbb
{S}^1}})$.
\begin{definition}  A knot (torus) $\hat \gamma~(\hat\lambda)$ in the manifold
$(S^2\times S^1,\eta_{_{\mathbb{S}^2\times\mathbb
{S}^1}})$ is said to be {\em trivial} if it is equivalent to the knot (torus)
$\hat\gamma_0~(\hat\lambda_0)$.
\end{definition}

\begin{proposition} \label{t.solid} Every $\eta^s_{_{S^2\times S^1}}$-essential
torus $\hat\lambda\subset(S^2\times S^1,
\eta^s_{_{S^2\times S^1}})$ bounds a solid torus in $S^2\times S^1$.
\end{proposition}

\begin{proposition} \label{kk}  A knot $\hat\gamma$ (torus $\hat\lambda$) in
the manifold $(S^2\times S^1,\eta^s_{_{S^2\times
S^1}})$ is trivial if and only if there is a tubular neighborhood
$N(\hat \gamma)~(N(\hat \lambda))$ of it in the manifold  $S^2\times S^1$
such that the manifold $(S^2\times S^1)\setminus N(\hat \gamma)~
((S^2\times S^1)\setminus N(\hat \lambda))$ is homeomorphic to the solid
torus (a pair of the solid tori).
\end{proposition}

\begin{figure}[h]
\centerline{\includegraphics[height=4cm]{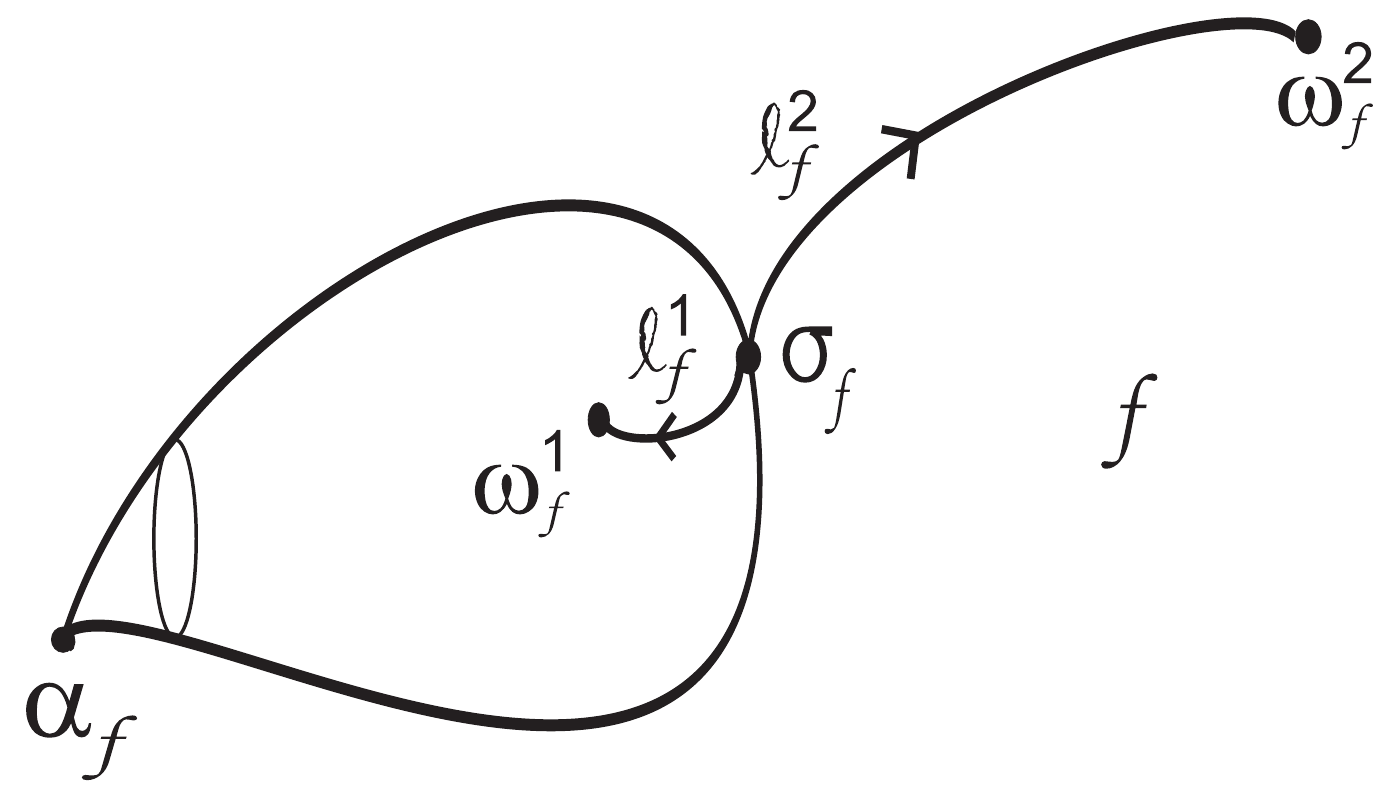}}\caption{{\small The phase portrait
of a diffeomorphism of the class $\mathcal P$}}
\label{r1}
\end{figure}

Denote by $\mathcal P$ the class of the Morse-Smale diffeomorphisms whose
non-wandering set consists of the source $\alpha_f$, the saddle $\sigma_f$
and the sinks $\omega_f^1 $, $\omega_f^2$. The phase portrait of a diffeomorphism
of the class $\mathcal P$ is shown in Figure \ref{r1}. Pixton has constructed  example
from class $\mathcal P$ mentioned above, so we call the class $\mathcal P$ the
{\it Pixton class}. We omit below the index $f$ in the notations of fixed
points.

A surprising fact is the existence of a countable set of non-conjugated diffeomorphisms
in the class $\mathcal P$. To understand this, we describe below knot topological
invariant suggested in \cite{BG}. Moreover, this invariant explains
existence in the class $\mathcal P$ of diffeomorphisms for which a saddle fixed point
possesses wildly embedded one-dimensional and two-dimensional separatrices.

Denote by $\ell_1$, $\ell_2$ the unstable 1-dimensional separatrices of the point $\sigma$.
It follows from Smale \cite{DDS} that the closure $cl(\ell_i)$ ($i=1,2$) is homeomorphic
to a simple compact arc which consists of the separatrix itself and two its extreme points:
$\sigma$ and the sink (see Proposition 2.3 in \cite{GMP}).
Moreover, the closures of the separatrices $\ell_1$ and $\ell_2$ contain different
sinks (see Corollary 2.2  in \cite{GMP}). To be definite, let $\omega_i$
belongs to $cl(\ell_i)$ (see Figure \ref{r1}).
For $i=1,2$ denote $V_i = W^s(\omega_i)\setminus\{\omega_i\}$. Denote by $\hat V_i =
V_i/f$ the corresponding orbit space and let $p_{_{i}}:V_i \to \hat V_i$ be the natural
projection that is the covering map inducing the epimorphism $\eta_{_{i}}:\pi_1(\hat V_i)
\to\mathbb Z$. As for the sink $\omega_i$ is concerned, the restriction
$f|_{V_i}$ is topologically conjugated to the diffeomorphism $a: \mathbb R^3\setminus \{O\}\to
{\mathbb R}^3\setminus \{O\}$, then the manifold $(\hat V_i,\eta_{_{i}})$ is equivalent to
the manifold $(S^2\times S^1,\eta^s_{_{S^2\times S^1}})$ and the set
$\hat\ell_i = p_{_{i}}(\ell_i)$ is the $\eta_{_{i}}$-essential knot in the manifold
$\hat V_i$ such that $\eta_{_{i}}(i_{_{\hat\ell_i*}}(\pi_1(\hat\ell_i)))=
\mathbb Z$ (see Theorem 2.3  in \cite{GMP}).

It was proved in \cite{BG} (Theorem 1) that at
 least one of the knots  $\hat\ell_1$, $\hat\ell_2$ is trivial (see also \cite{GMP},
Proposition 4.3).
To be definite we assume below the knot  $\hat\ell_1$ be trivial.

Next result was proved in \cite{BG} (Theorem 3) (see also \cite{GMP},
Theorem 4.3).

\begin{proposition} Diffeomorphisms $f,f'\in \mathcal P$ are topologically conjugated
if and only if the knots $\hat\ell_2(f)$ and $\hat\ell_2(f')$ are equivalent.
\label{iffknot}
\end{proposition}

Therefore the equivalence class of the knot $\hat\ell_2(f)$ is a complete topological
invariant for diffeomorphisms from the Pixton class. Moreover, the following realization
theorem holds (see \cite{BG} Theorem 2 and \cite{GMP} Theorem 4.4).

\begin{proposition}\label{Exist} For every knot $\hat\ell\subset(S^2\times
S^1, \eta^s_{_{S^2\times S^1}})$ such that
$\eta^s_{_{S^2\times S^1}}(i_{_{\hat\ell*}}(\pi_1(\hat\ell)))=
\mathbb Z$ there is a diffeomorphism  $f: S^3\to S^3$ from the class
$\mathcal P$ such that the knots $\hat\ell$ and $\hat\ell^2(f)$ are equivalent.
\label{t.three}
\end{proposition}

Masur constructed an example of an essential and nontrivial knot embedded
to $S^2\times S^1$\cite{Mazur}. According  to Proposition \ref{Exist}, there exists
a diffeomorphism $f$ of the Pixton class such that exactly one unstable one-dimensional
separatix and stable two-dimensional separatix of the saddle point $\sigma$ are
wildly embedded.

\begin{figure}[h]
\centerline{\includegraphics[height=6cm]{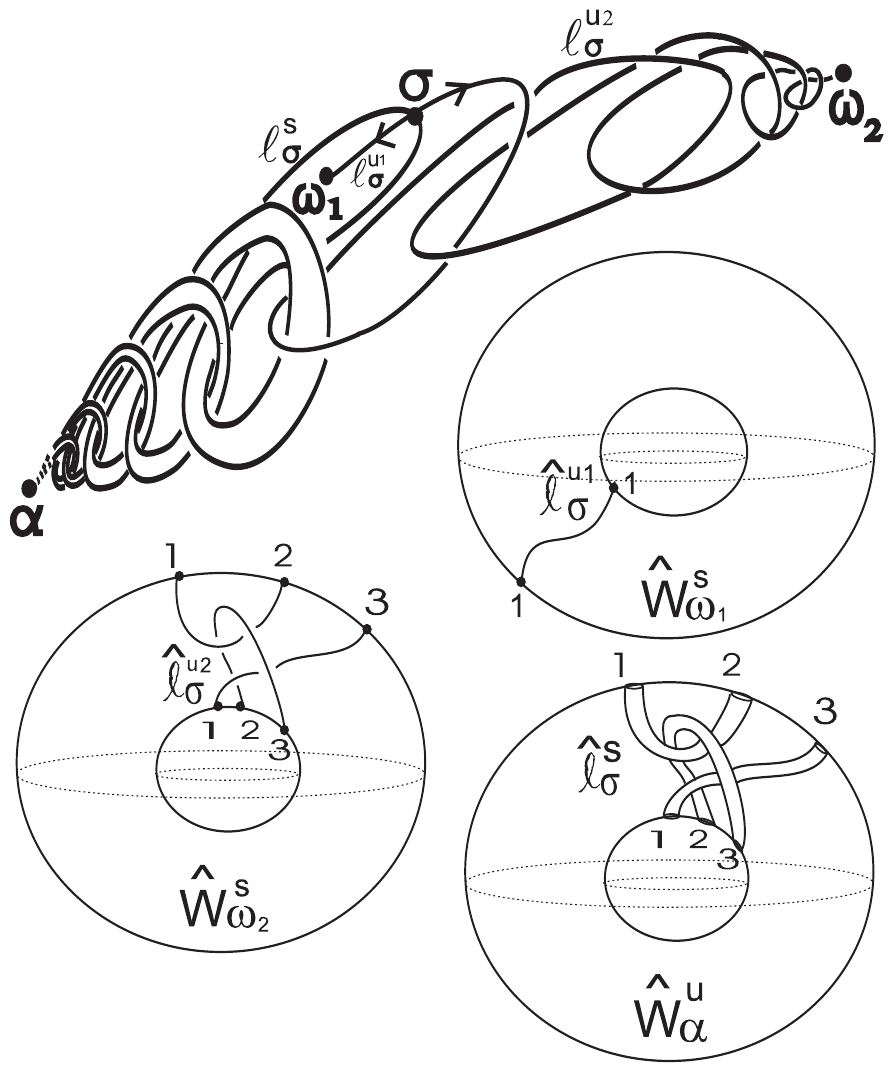}}\caption{{\small The phase portrait
of a diffeomorphism of the class $\mathcal P$ and projection of saddle separatrices
in factor spaces}}
\label{r2}
\end{figure}

On Fig.~\ref{r2} it is shown the Masur's knot $\hat l^u_\sigma$  which appears in
factor-space $\hat W^s(\omega^2)$ and essential torus $\hat l^s_\sigma$ embedded to
$\hat W^u(\alpha)$ which is tubular neighborhood of the a Masur knot.

\subsection{ Diffeomorphisms with wildly embedded frames}

\begin{definition}
For $k\in\mathbb N$ a $k$-frame $F_k$ in $\mathbb R^n$ at the point $p$ is an union
of $k$ simple curves $A_1, \ldots, A_k$, $F_k=\bigcup\limits_{i=1}^k A_i$,
with a single common point $p$ such that $p$ is the
boundary point of each $A_k$, $k\ge 1$ and $A_i\cap A_j= p,~i\neq j$.
\end{definition}

\begin{definition} $ $
\begin{itemize}
\item The $k$-frame $F_k=\bigcup\limits_{i=1}^kA_i$ is said to be {\em standard}
if each arc $A_i$ lies in the plain $Ox_1x_2$ and it is defined $\varphi=\frac{2\pi(i-1)}
{k}$ where $\rho,~\varphi$ are the polar coordinates in the plain $Ox_1x_2$.
\item A $k$-frame $F_k=\bigcup\limits_{i=1}^kA_i$ is said to be {\em
tame}, if there is a homeomorphism $\varphi:\mathbb R^n\to\mathbb R^n$ such that
$\varphi(F_k)$ is standard. Otherwise, the frame $F_k$ is said to be {\em wild}.
\item A $k$-frame $F_k=\bigcup\limits_{i=1}^k A_i$ is said to be {\em mildly wild},
if the frame $F_k\setminus (A_i\setminus p)$ is tame for every $i\in\{1,\dots,k\}$.
\end{itemize}
\end{definition}

\begin{figure}[h]
\centerline{\includegraphics[height=6cm]{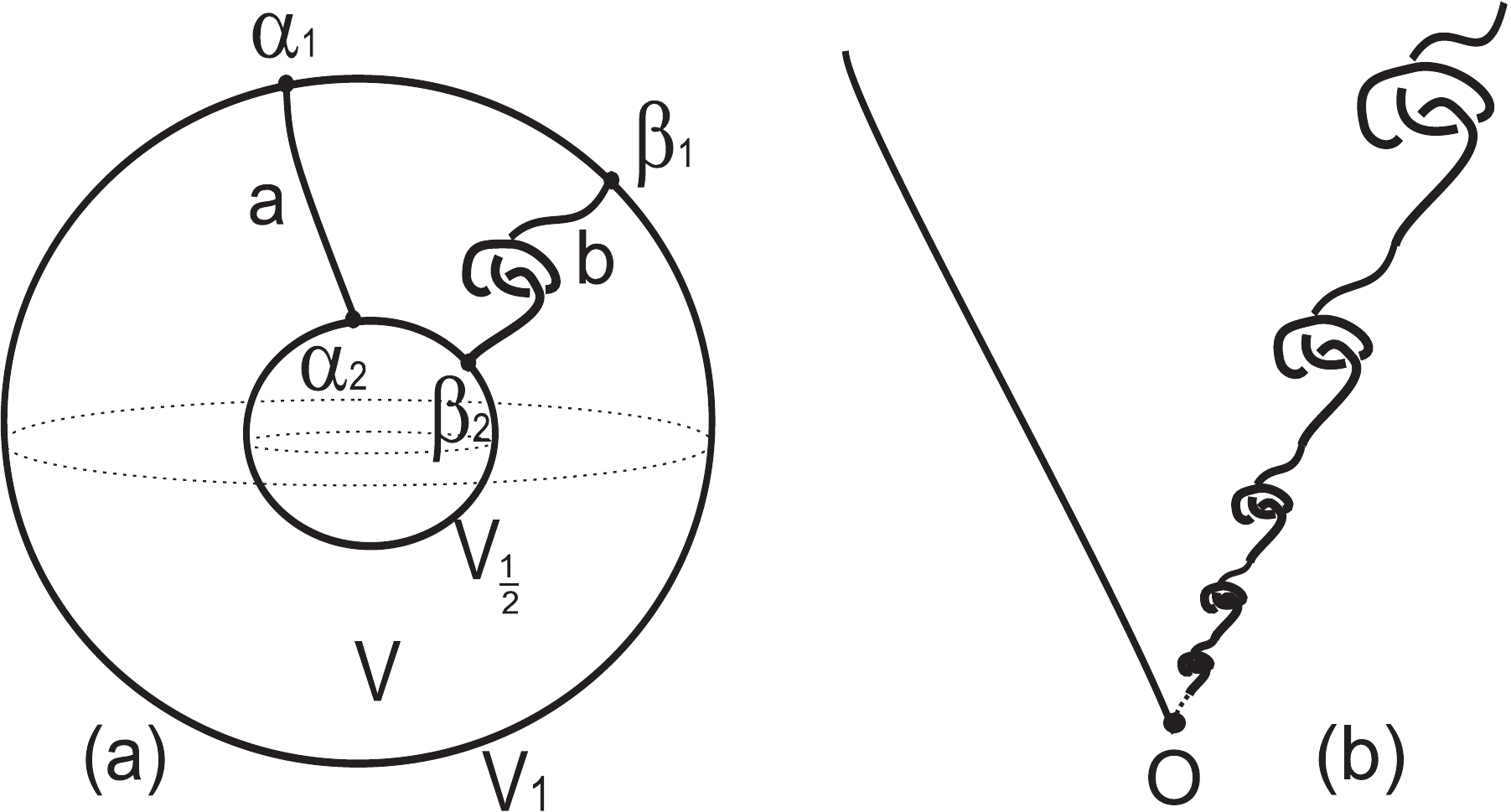}}\caption{{\small
A construction of a wild 2-frame in $\mathbb R^3$}}
\label{r5}
\end{figure}

\begin{figure}[h]
\centerline{\includegraphics[height=6cm]{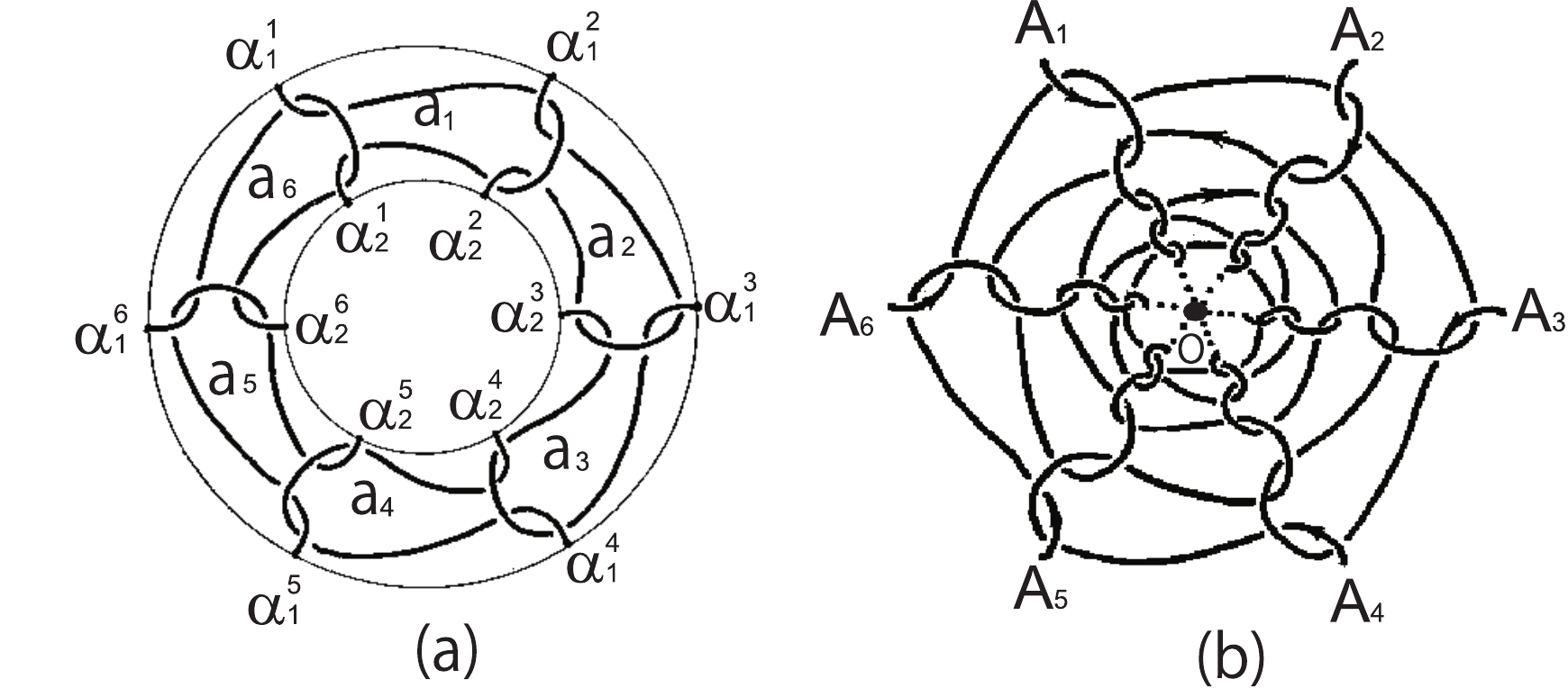}}\caption{{\small Debrunner-Fox's
example}}
\label{r6}
\end{figure}

One can easily construct a wild $k$-frame, if one assumes the arc $A_1$ be the wild arc
$\tilde\ell$ of Artin-Fox's example \cite{AF}. But the fact that each arc $A_i$ is tame
does not mean that the frame $F_k$ is tame. Figure \ref{r5} (b) shows
an example of the wild 2-frame.  Similarly to the Artin-Fox's example this frame is
constructed using the arcs $a,~b$ shown in Figure \ref{r6} (a).
The boundary points $\alpha_1,~\alpha_2,~\beta_1,~\beta_2$ of the respective arcs
$\alpha,~\beta$  are glued by $\phi(\alpha_1)=\alpha_2$, $\phi(\beta_1)=\beta_2$ and
$A_1=\bigcup\limits_{k\in\mathbb Z}\phi^k(a)\cup O$,
$A_2=\bigcup\limits_{k\in\mathbb Z}\phi^k(b)\cup O$, $F_2=A_1\cup A_2$.
From the Statement \ref{tame-arc} it follows that both $A_1,~A_2$ are tame.
Debrunner and Fox \cite{DF} presented the construction of a mildly
wild $k$-frame for every $k>1$. Figure \ref{r6} shows this construction
for $k=6$.  This 6-frame is constructed using the arcs $a_1,\dots,a_6$ shown in
Figure \ref{r6} (a). The boundary points $\alpha_1^i,~\alpha_2^i,$ of the
arc $a_i,~i\in\{1,\dots,6\}$ are glued by $\phi(\alpha_1^i)=\alpha_2^i$ and
$A_i=\bigcup\limits_{k\in\mathbb Z}\phi^k(a_i)\cup O$, $F_6=\bigcup\limits_{i=1}^6A_i$.

Let $f$ be a Morse-Smale diffeomorphism and suppose that there is a sink
$\omega\in NW(f)$ and the set $L_\omega$  of all unstable one-dimensional different
separatrices $\ell_1, \dots, \ell_k$ of saddles $\sigma_1, \dots, \sigma_r$, $k, r
\in \mathbb N$, $k\leq r$ ($\sigma_i$ may coincide
with $\sigma_j$) such that for any $j$ closure of separatrix $\ell_j$ consists of
exactly two points: $\omega$  and saddle point $\sigma$ for which $\ell_j$ is separatrix.
Since $W^s_\omega$ is homeomorphic to $\mathbb R^3$ and since the set
$L_\omega\cup\omega$ is the union of the simple arcs with the unique common point
$\omega$ belonging to each arc, analogously to a frame of arcs in  $\mathbb R^3$ we
call $L_\omega\cup\omega$ the {\it frame of 1-dimensional unstable separatrices}.

\begin{definition} A frame of separatrices $L_\omega\cup\omega$ is {\em tame} if there
is a homeomorphism $\psi_\omega: W^{s}_\omega \to \mathbb R^3$ such that
$\psi_\omega (L_\omega\cup\omega)$ is the standard frame of arcs in $\mathbb R^3$.
Otherwise the frame of separatrices is called {\em wild}.
\end{definition}

\begin{figure}[h]
\centerline{\includegraphics[height=6cm]{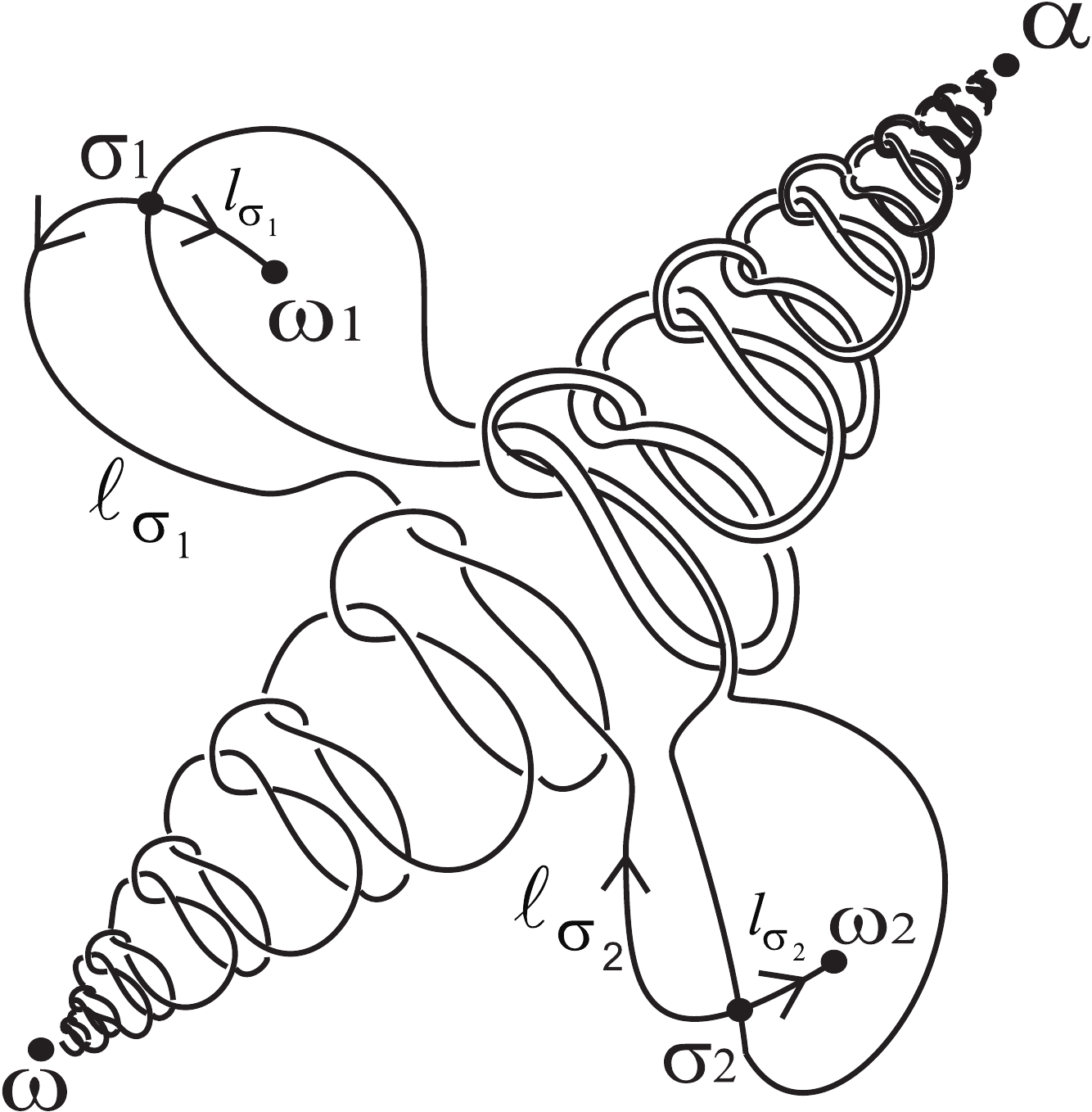}}\caption{{\small Phase portrait of a
Morse-Smale diffeomorphism on $\mathbb S^3$ with the mildly wild frame of separatrices}}
\label{r4}
\end{figure}

In \cite{Pochinka2009} by the method, similar to that described in subsection
\ref{Pixton}, a Morse-Smale diffeomorphism was constructed having a mildly wild frame
of one-dimensional separatrices (see fig. \ref{r4}).

\section{Periodic vector field on $S^3$ with wildly embedded \\separatrix set}

Now we present a periodic vector field on $S^3$ with a wild embedding of
a 2-dimensional unstable separatrix manifold and 3-dimensional stable separatrix
manifold for the saddle IC with exponential dichotomy on $\R$ of the type $(3,2).$
Also we present another periodic vector field on $S^3$ that has a mildly wild frame
of 2-dimensional separatrix manifolds.

We start with some diffeomorphism $f$ of the Pixton class on $S^3$ that
has one hyperbolic source $\alpha$, one saddle $\sigma$ of the type $(2,1)$ (2-dimensional
stable manifold and 1-dimensional unstable one) and two hyperbolic sinks $\omega_1, \omega_2$.
Stable 2-dimensional manifold of $\sigma$ contains in its closure the point $\alpha$, that is,
all orbits of $f$ with initial points on $W^s(\sigma)$, except $\sigma$ itself, have the only
$\alpha$-limit point $\alpha$ and the $\omega$-limit point $\sigma$. The closure
of $W^s(\sigma)$ is a topologically embedded sphere $\Sigma$ in $S^3$,
being the boundary of two open 3-balls $D_1$, $D_2$ in $S^3$. The fixed point $\omega_1$
(sink) lies inside of the ball $D_1$, another sink $\omega_2$ lies inside another ball $D_2$.
We suppose that a 1-dimensional separatrix of $\sigma$ which enters to $D_2$ is wildly
embedded. This implies stable manifold $W^s(\sigma)$ be also wildly embedded (see figure \ref{p}).

Now consider the suspension over $f$. Since $f$ is diffeotopic to
$id_{S^3},$ the manifold $M_f$ is topologically a direct product $S^3\times S^1$,
moreover, this direct product structure can be chosen by means of a some diffeomorpism
(see above). We fix this product structure and consider henceforth the suspension
as the standard $S^3\times S^1$. Thus, the suspension flow in $S^3\times S^1$ has
one totally unstable periodic orbit, one saddle periodic orbit of the type $(3,2)$
and two totally stable periodic orbits, all of them are hyperbolic periodic orbits.
The projection of any of these periodic orbits onto the base $S^1$ is 1-1 correspondence.

Now recall that the suspension flow is Morse-Smale one. All its periodic orbits are
hyperbolic and any other orbit tends to some of four periodic orbits as $t\to \pm\infty$.
This implies, by the construction, all four periodic ICs of the non-autonomous vector field
on $S^3$, to possess an exponential dichotomy on $\R$.
Types of an exponential dichotomy are different: two stable periodic orbits give
rise to two completely stable periodic ICs, their type of an exponential dichotomy
is $(4,1)$ (four is the dimension of their stable manifolds),
the saddle periodic orbit gives rise to the saddle periodic IC with the dichotomy
on $\R$ of the type $(3,2),$ and the completely unstable periodic orbit
gives rise to the IC with the dichotomy of the type $(1,4)$. All other ICs tend
to these four ICs and hence they possess an exponential dichotomy on $\R_-$ and $\R_+$
separately depending on which of four ICs they approach to.

Now remind that diffeomorphism $f$ has a smooth curve being the unstable
separatrix $W^u(\sigma)$ for the saddle point $\sigma$. For the suspended flow in
$S^3\times S^1$ we get the two-dimensional smooth unstable submanifold $W^u(\gamma_\sigma)$
of the saddle periodic orbit $\gamma_\sigma$. The manifold $W^u(\gamma_\sigma)$ is a direct
product of $W^u(\sigma)\times S^1$, this follows from the suspension construction.
If one of two unstable separatrices of $\sigma$ is wildly embedded into $S^3$ (see above),
then one of connected component of intersection  $S^3_\tau \cap W^u(\gamma_\sigma),$
$\tau \in S^1,$ (denote it as $\Sigma_\tau$) is a wildly embedded curve in $S^3_\tau$.
Suppose, to be definite, that stable sink $\omega_2$ be the $\omega$-limit set for
all orbits of the wildly embedded separatix of $\sigma$. We will say that the related
component of $W^u(\sigma)\times S^1\setminus \gamma_\sigma$ is wildly
embedded in $S^3\times S^1$. The characterization of this wild embedding is the
following. Choose any smooth 3-disk $D$ being transversal to a point on the periodic
orbit $\gamma_{\omega_2}$. Then the wildly embedded component intersects
this disk along a smooth ray with the extreme point $D\cap\gamma_{\omega_2}.$

\begin{lemma}
If flows $f^t, f^{'t} \in  \tilde{\mathcal P}$ are topologically
equivalent and $f^t$ possesses a wildly embedded connected component of the
set $W^u(\gamma_S)\setminus\gamma_S$, then the same holds true for
the flow $f^{'t}$.
\end{lemma}

Thus, we have proved the assertion
\begin{theorem}
There exists a smooth 1-periodic vector field $v$ on $S^3$ such that $v$ is
gradient-like one with only four 1-periodic ICs possessing exponential dichotomies
on $\R$: completely unstable (of the type (1,4)), saddle one of the type $(3,2),$
two completely stable ones and its saddle periodic IC has wildly embedded two-dimensional
and three-dimensional separatrix sets.
\end{theorem}

In the same way, starting with a diffeomorphism of $S^3$ having a mildly
wild frame of separatrices described above, we get a 1-periodic gradient-like vector
field on $S^3$ such that it has one completely unstable IC, one completely
stable IC and $n\ge 2$ saddle periodic ICs with an exponential dichotomy
of the type (3,2) whose $n$ two-dimensional unstable separatrices form a mildly
wild frame along with their $n$ three-dimensional stable separatrices which also form
a mildly wild frame similar to that plotted in Fig.\ref{r4}.

We formulate some assertion concerning non-autonomous vector fields we
have constructed.
\begin{theorem}\label{rough}
Any sufficiently small uniform perturbation of such vector field $v$ gives a
non-autonomous vector field $v'$ that is uniformly equivalent to the initial
one, i.e. there exists an equimorphism $h$ of the extended phase manifold $S^3\times
\R$ which transforms the foliation $\mathcal L_v$ to that of $\mathcal L_{v'}.$
\end{theorem}
Due to the very simple structure of $\mathcal L_v$ the proof is almost
evident, nevertheless, because of some technicalities, it will be performed elsewhere.

\section{Perturbations}

In this section we perturb periodic vector fields constructed in
the preceding section in such a way that its uniform structure stays the
same, but, in dependence on the perturbation chosen, the perturbed non-autonomous
vector field would be almost periodic or even be nonrecurrent in time.

By the theorem \ref{rough}, $v$ is structurally stable w.r.t. small
uniform perturbations of the form $v+\varepsilon v_1$ given by a bounded
uniformly continuous map $v_1:\R \to V^r(S^3)$ into the Banach space $V^r(S^3)$.
In particular, such a perturbation can be
chosen being an almost periodic map. In this case, since any of four periodic ICs
possess an exponential dichotomy on $\R$ (of different types) the perturbed almost
periodic vector field will have in a small enough uniform neighborhood of each
periodic IC an almost periodic IC with the same type of exponential
dichotomy. Moreover, the perturbed vector field will be also gradient-like
one and any of its IC will possess an exponential dichotomy on $\R_{\pm}$ and
will tend to some of four almost periodic ICs.

The periodic vector field, we constructed above, has a very simple structure of
its foliation into ICs. Namely, it has a unique totally unstable IC
$\gamma_\alpha$ possessing an exponential dichotomy on $\R$ of the type (4,1), one
saddle IC $\gamma_\sigma$ possessing an exponential dichotomy on $\R$ of the type
(3,2) and two totally stable ICs $\gamma_{\omega_1}$ and
$\gamma_{\omega_2}$ possessing both an exponential dichotomy on $\R$ of the type
(4,1). All other ICs tend to one of these specific ICs as $t\to \pm
\infty$. One important thing exists. Let us choose some sufficiently thin
uniform neighborhoods $U_j,$ $j=1-4$, of all specific ICs. One can choose
these neighborhoods in such a way that the passage time within
$M\times \R\subset \cup_jU_j$ would be uniformly bounded from above and
below. This allow us to prove that the non-autonomous vector field is
structurally stable with respect to a small enough uniformly bounded
perturbations. This means that there is an equimorphism
$\Phi:M\times \R \to M\times \R$
such that it transforms a foliation into ICs of the periodic
vector filed into the foliation of the perturbed vector field. In
particular, $\Phi$ preserves the properties of a wild embedding or mildly
wild separatrix frame. If $v_1$ is nonrecurrent but uniformly continuous
bounded map, then the perturbed NVF will have the same uniform structure
but its four specific ICs will lie in thin uniform neighborhoods of those
for the constructed periodic NVF. The same holds true for the periodic NVF
with the mildly wild $k$-frame of separatrices.

\section{Addendum: Elements of the uniform topology}
\label{Addendum}

For the reader convenience, we present here some notions from the uniform topology.
Recall some basic definitions of the theory of uniform spaces (see details in
\cite{Kelley}). A set $X$ is called an {\em uniform space}, if on $X\times
X$ is defined a collection $\mathcal U$ of its subsets satisfying the following
conditions (if so, $\mathcal U$ is called {\em the uniformity})
\begin{enumerate}
\item each element of $\mathcal U$ contains diagonal $\Delta = \cup_{x\in X}\{(x,x)\};$
\item if $U\in \mathcal U$, then $U^{-1}\in \mathcal U$, where $U^{-1}$ is
the set of all pairs $(y,x)$ for which $(x,y)\in U;$
\item for any $U\in \mathcal U$ some $V\in \mathcal U$ exists such that $V\circ V \subset
\mathcal U$, here $V\circ V$ denotes the composition: $(x,z)\in V\circ V$ if there is
$y\in X$ such that $(x,y)\in V$ and $(y,z)\in V;$
\item if $U,V \in \mathcal U,$ then $U\cap V \in \mathcal U;$
\item if $U \in \mathcal U$ and $U\subset V \subset X\times X,$ then $V\in \mathcal U.$
\end{enumerate}
If $X$ is a metric space with metrics $d$, then 1) corresponds to the property
$d(x,x)=0,$ 2) corresponds to the symmetry of $d$: $d(x,y)=d(y,x).$ The
property 3) is of the type of the triangle inequality: for any ball of the
radius $r$ a ball of the radius $r/2$ should exist. The third and fifth
conditions are similar to the axioms of neighborhoods near a point for the
topology that is defined by the uniformity.

The uniformity on a given set $X$ can be defined by many ways providing
different uniform spaces. This was used above where on the set $M\times \R$
different uniform structures were defined. When $(X,\mathcal U),$ $(Y, \mathcal V)$ are
two uniform spaces, then a notion of a uniformly continuous map $h: X\to
Y$ is defined. Namely, a map $h: X\to Y$ is uniformly continuous w.r.t.
$\mathcal U, \mathcal V$, if for any $V\in \mathcal V$ the set
$\{(x,y)| (h(x),h(y))\in V\}$ belongs to $\mathcal
U$. When $h: X\to Y$ is one-to-one and both $h, h^{-1}$ are uniformly
continuous, then $h$ is called to be an equimorphism. In
this case uniform spaces $(X,\mathcal U),$ $(Y, \mathcal V)$ are called
uniformly equivalent or equimorphic ones.

An uniformity on the set $X$ making it the uniform space $(X,\mathcal
U)$ generates the definite topology on $X$ making it a topological space.
This space can possess various topological properties. Conversely, each
regular\footnote{A topological space is regular, iff for every
its point $x$ and any its neighborhood $U$ there is a closed neighborhood
$V$ of $x$ such that $V\subset U$.} topology $\mathcal T$ on $X$ is an uniform
topology which corresponds to some uniformity, but such uniformity is, in general,
not unique. But if the topological space is compact and regular one, then there is
an unique uniformity generating the topology $\mathcal T$.

\section{Acknowledgement}

This research was supported by the Laboratory of Dynamical Systems and Applications
of NRU HSE, of the Ministry of Science and Higher Education of the RF grant agreement
075-15-2019-1931. L.M. Lerman also acknowledges some support from the Research and Educational
Center ``Mathematics for Future Technologies'' (agreement 075-02-2021-1394).

\section{Conflict of Interest}
The authors declare they have not conflict of interest.

\end{document}